\documentclass[12pt]{amsart}

\usepackage{amsmath,amssymb,array,enumerate,url}

\usepackage[canadien]{babel}
\usepackage[utf8]{inputenc}
\usepackage[T1]{fontenc}	

\newtheorem{thm}{Th\'eor\`eme}

\newcommand {\abs}[1]{\left| #1\right|}
\newcommand{\bZ}{\mathbb{Z}}
\newcommand{\bQ}{\mathbb{Q}}

\newcommand{\bR}{\mathbb{R}}
\newcommand{\cO}{{\mathcal{O}}}
\newcommand{\cS}{\mathcal{S}}
\newcommand{\et}{\quad\mbox{et}\quad}

\newcommand{\GL}{\mathrm{GL}}
\newcommand{\Kbar}{\overline{K}}
\newcommand{\pgcd}{\mathrm{pgcd}}
\newcommand{\tX}{{\tilde{X}}}
\newcommand{\ux}{\boldsymbol{x}}

\DeclareMathOperator{\trace}{trace}

\begin{document}

\baselineskip=14.6pt

\title[Construction de nombres extr\'emaux]{Construction de nombres extr\'emaux\\ pour le probl\`eme de l'approximation simultan\'ee\\ d'un nombre et de son carr\'e}

\author{Samuel Pilon \et Damien Roy}


\maketitle

\begin{center}
    \begin{minipage}{14 cm} {\bf R\'esum\'e.}
        Nous consid\'erons le probl\`eme de l'approximation simultan\'ee d'un nombre et de son carr\'e dans un cadre g\'en\'eral qui englobe aussi bien les corps quadratiques imaginaires que les corps de fonctions rationnelles en une variable. Dans ce contexte, nous construisons de nouveaux nombres extr\'emaux.
    \end{minipage}
\end{center}

\smallskip
\begin{center}
    \begin{minipage}{14 cm} {\bf Abstract.}
{\sl Construction of extremal numbers for the problem of simultaneous approximation to a number and to its square.}
        We consider the problem of simultaneous approximation to a number and to its square in a general framework that encompasses imaginary quadratic number fields and fields of rational functions in one variable. In this context, we construct new extremal numbers.
    \end{minipage}
\end{center}

\bigskip
Le probl\`eme de l'approximation simultan\'ee d'un nombre et de son carr\'e a \'et\'e initialement consid\'er\'e par Davenport et Schmidt comme outil pour l'approximation des nombres r\'eels par les entiers alg\'ebriques cubiques. En 2013, Bel a \'etendu le cadre de ce probl\`eme à un corps de nombres $K$ quelconque.

Pour son r\'esultat principal, Bel choisit, pour chaque place $v$ de $K$, la valeur absolue $\abs{\cdot}_v$ qui \'etend la valeur absolue usuelle de $\bQ$ si $v$ est archim\'edienne, ou sa valeur absolue $p$-adique usuelle si $v$ est au-dessus d'un nombre premier $p$ (on a alors $|p|_v=p^{-1}$). Il fixe aussi une place $w$ de $K$ et d\'esigne par $K_w$ le compl\'et\'e de $K$ pour cette valeur absolue et par $\cO_{K,w}$ l'ensemble des \'el\'ements $x$ de $K$ tels que $|x|_v\leq 1$ pour toute place $v$ de $K$ distincte de $w$. En notant par $\gamma=(1+\sqrt{5})/2$ le nombre d'or, le th\'eor\`eme 1.1 de \cite{Be} s'\'enonce ainsi.

\begin{thm}[Bel]\label{Bel}
Soient $K$ et $w$ comme ci-dessus.
    \begin{enumerate}[(a)]
        \item Pour chaque $\xi\in K_w$ tel que $1,\xi,\xi^2$ soient lin\'eairement ind\'ependants sur $K$, il existe une constante $c>0$ telle que les conditions
        \[
         |x_0|_w\leq X
         \et
         \max\{|x_0\xi-x_1|_w, |x_0\xi^2-x_2|_w\}\le c X^{-1/\gamma}
        \]
        n'admettent pas de solution non nulle $(x_0,x_1,x_2)\in\cO^3_{K,w}$ pour des nombres r\'eels $X$ arbitrairement grands.
        \item Pour chaque $\epsilon>0$, il existe un nombre $\xi\in K_w$ qui est transcendant sur $K$ tel que les conditions
        \[
         |x_0|_w \le X
         \et
         \max\{ |x_0\xi-x_1|_w, |x_0\xi^2-x_2|_w \}\le  X^{-1/\gamma+\epsilon}
        \]
        admettent une solution non nulle $(x_0,x_1,x_2)\in\cO^3_{K,w}$
        pour tout nombre r\'eel $X$ assez grand.
    \end{enumerate}
\end{thm}

Dans les travaux ant\'erieurs à Bel, on n'avait consid\'er\'e que le corps $K=\bQ$. Dans ce cas, si $w=\infty$ on a $K_w=\bR$ et $\cO_{K,w}=\bZ$. La partie (a) du th\'eor\`eme reprend alors le contenu du th\'eor\`eme 1a de \cite{DaSc}, tandis que la partie (b) est d\'emontr\'ee sous une forme plus forte au th\'eor\`eme 1.1 de \cite{Ro1}, avec la borne $X^{-1/\gamma+\epsilon}$ remplac\'ee par $c'X^{-1/\gamma}$ pour un $c'>0$, les nombres r\'eels $\xi$ qui r\'ealisent cette borne \'etant dits \emph{extr\'emaux}.

Si on prend pour $w$ la place de $\bQ$ associ\'ee à un nombre premier $p$, on a $K_w=\bQ_p$ et $\cO_{K,w}=\bZ\left[1/p\right]\cap[-1,1]$. La partie (a) est alors \'equivalente au cas $n=3$ du th\'eor\`eme 2 de \cite{Te}, tandis que la partie (b) est prouv\'ee ind\'ependamment par Zelo au th\'eor\`eme 1.3.4 de \cite{Ze} et par Bugeaud dans \cite{Bu}. Dans ce cas, on ne connait pas de nombres extr\'emaux et on doute de leur existence.

Le but de cet article est de pr\'esenter de nouvelles situations pour lesquelles il existe des nombres extr\'emaux, en cons\'equence du r\'esultat suivant.

\begin{thm}
Soient $A$ un anneau int\`egre et $K$ son corps des fractions.  On suppose qu'il existe une valeur absolue non triviale $\abs{\cdot}$ sur $K$ telle que $\abs{a}\ge 1$ pour tout $a\in A\setminus\{0\}$, et on  note $\Kbar$ le compl\'et\'e de $K$ pour cette valeur absolue.
    \begin{enumerate}[(a)]
    \item Supposons que $A$ soit un anneau \`a factorisation unique.  Alors, pour tout $\xi\in\Kbar$ tel que 1, $\xi$ et $\xi^2$ sont lin\'eairement ind\'ependants sur $K$, il existe une constante $c_1>0$ et des nombres r\'eels $X$ arbitrairement grands tels que les in\'egalit\'es
        \begin{equation}
         |x_0|\le X
         \et
         L(\ux):=\max\{ |x_0\xi-x_1|,\, |x_0\xi^2-x_2|\}
         \le c_1X^{-1/\gamma}
         \label{ineq1}
        \end{equation}
        ne poss\`edent aucune solution non nulle $\ux=(x_0,x_1,x_2)\in A^3$.
    \item Il existe $\xi\in\Kbar$ tel que 1, $\xi$ et $\xi^2$ sont lin\'eairement ind\'ependants sur $K$ et une constante $c_2>0$ telle que pour tout $X$ assez grand les in\'egalit\'es
        \begin{equation}
          |x_0|\le X
          \et
          L(\ux):=\max\{ |x_0\xi-x_1|,\, |x_0\xi^2-x_2|\}
          \le c_2X^{-1/\gamma}
          \label{ineq2}
        \end{equation}
        poss\`edent une solution non nulle $\ux=(x_0,x_1,x_2)\in A^3$.
    \end{enumerate}
\end{thm}

Par exemple, on peut prendre pour $A$ l'anneau des polyn\^omes $F[u]$ en une ind\'etermin\'ee sur un corps quelconque $F$.  En munissant cet anneau de la valeur absolue $|p|=e^{\deg(p)}$ et en \'etendant cette derni\`ere au corps $K=F(u)$, on obtient pour compl\'et\'e $\Kbar=F((1/u))$. Alors le th\'eor\`eme montre l'existence de s\'eries $\xi\in F((1/u))$ avec une propri\'et\'e d'approximation extr\'emale.

Lorsque $K$ est un corps de nombres, la condition $|a|\ge 1$ pour tout $a\in A$ implique que $K$ ne poss\`ede qu'une place archim\'edienne $w$, que $\abs{\cdot}$ est une valeur absolue associ\'ee \`a cette place et que $A$ est contenu dans l'anneau des entiers $\cO_K=\cO_{K,w}$ de $K$.  Alors, selon la partie (a) du th\'eor\`eme de Bel cit\'e pr\'ec\'edemment, le th\'eor\`eme 2 (a) demeure vrai sans supposer que $A$ soit \`a factorisation unique.  Donc le th\'eor\`eme 2 (b) montre l'existence de nombres complexes extr\'emaux sur l'anneau des entiers d'un corps quadratique imaginaire donn\'e, par exemple sur $\bZ[i]$ ou sur $\bZ[\sqrt{-5}]$.


\subsection*{D\'emonstration de la partie (a).}
On reprend la d\'emonstration de \cite[Theorem 1a]{DaSc} en supposant que les in\'egalit\'es \eqref{ineq1} admettent une solution non-nulle $\ux=(x_0,x_1,x_2)$ dans $A^3$ pour tout $X$ assez grand, disons pour tout $X\ge X_0$.  En vertu des hypoth\`eses sur $A$, on peut se restreindre aux solutions \emph{primitives}, c'est-\`a-dire celles pour lesquelles $\pgcd(x_0,x_1,x_2)=1$.  Le but est de montrer que $c_1>c(\xi)>0$.  Cela permet de supposer que $c_1<6^{-1/2}$ et que $X_0\ge 1$.

Pour chaque $X\ge X_0$, on note $\cS(X)$ l'ensemble non vide des solutions primitives de \eqref{ineq1} dans $A^3$, et on pose $\ell(X)=\inf\{L(\ux)\,;\,\ux\in\cS(X)\}$.  Comme $L(\ux)<(6X)^{-1/2}$ pour tout $\ux\in \cS(X)$, les estimations de \cite[Lemma 4]{DaSc} montrent que le d\'eterminant de n'importe quels trois points de $\cS(X)$ est en valeur absolue inf\'erieur \`a $1$, donc nul.  Ainsi, $\cS(X)$ est contenu dans un sous-$K$-espace vectoriel $V$ de $K^3$ de dimension au plus $2$ et l'argument de \cite[Lemma 5]{DaSc} permet de conclure que $\ell(X)>0$.  On peut donc construire r\'ecursivement une suite de points primitifs $(\ux_i)_{i\ge 1}$ dans $A^3$ en choisissant, pour chaque $i\ge 1$,
\[
 \ux_i\in\cS(\tX_{i-1})
 \quad \text{tel que}\quad
 L(\ux_i) < 2\ell(\tX_{i-1})
\]
o\`u $\tX_{i-1}=X_0$ si $i=1$ et $\tX_{i-1}=\big(2c_1/L(\ux_{i-1})\big)^\gamma$ si $i>1$.  Alors, en posant $X_i = |x_{i,0}|$ et $L_i = L(\ux_i)$ pour tout $i\ge 1$, on obtient
\begin{equation}
\label{ineq3}
 L_i\le 2c_1X_{i+1}^{-1/\gamma}, \quad L_{i+1}\le L_i/2 \et X_i\le X_{i+1}
 \quad
 (i\ge 1).
\end{equation}
En effet, comme $\ux_{i+1}\in\cS(\tX_i)$, on a $X_{i+1}\le \tX_i$, donc $L_i=2c_1\tX_i^{-1/\gamma}\le 2c_1X_{i+1}^{-1/\gamma}$.  Comme $\ux_{i+1}\in\cS(\tX_i)$, on a aussi $L_{i+1}\le c_1\tX_i^{-1/\gamma}=L_i/2$.  Enfin, puisque $L_i<2\ell(\tX_{i-1})$, cette derni\`ere in\'egalit\'e implique $L_{i+1}<\ell(\tX_{i-1})$, donc $\ux_{i+1}\notin\cS(\tX_{i-1})$, et par suite $X_{i+1}>\tX_{i-1}\ge X_i$.

Les in\'egalit\'es \eqref{ineq3} impliquent que $\lim_{i\to\infty}X_i=\infty$ car si la suite $(X_i)_{i\ge 1}$ \'etait born\'ee sup\'erieurement par un nombre r\'eel $X$, on aurait $\ell(X)\le L_i$ pour tout $i\ge 1$, donc $\ell(X)=0$, ce qui est impossible.  Comme $L_i\neq L_{i+1}$, on observe aussi que les points $\ux_i$ et $\ux_{i+1}$ sont lin\'eairement ind\'ependants sur $K$ pour tout $i\ge 1$.  Partant de l\`a, on conclut en reprenant, avec des changements mineurs, le reste de l'argument de Davenport et Schmidt dans \cite[\S 3]{DaSc}.  Par exemple si, pour un entier $i\ge 1$, on a $x_{i,0}x_{i,2}-x_{i,1}^2=0$, alors $\ux_i=a(m^2,mn,n^2)$ pour une unit\'e $a$ de $A$ et un point non nul $(m,n)$ de $A^2$. Comme $|a|=1$, on en d\'eduit que $|m|=X_i^{1/2}$ et que $|m\xi-n|\le X_i^{-1/2}L_i$.

\subsection*{D\'emonstration de la partie (b).}
Soit $\rho>1$ et soit $(a_j)_{j\ge 1}$ une suite d'\'el\'ements de $A$ avec $|a_j|\ge 1+\rho$ pour tout $j\ge 1$.  On peut donner un sens \`a la fraction continue $\xi=[0,a_1,a_2,a_3,\dots]$ comme limite dans $\Kbar$ des quotients $p_j/q_j$ o\`u $(p_j)_{j\ge 0}$ et $(q_j)_{j\ge 0}$ sont les suites d'\'el\'ements de $A$ caract\'eris\'ees par
\[
 \begin{pmatrix} q_j &q_{j-1}\\ p_j &p_{j-1}\end{pmatrix}
 =
 \begin{pmatrix} a_1 &1\\ 1 &0\end{pmatrix}
 \begin{pmatrix} a_2 &1\\ 1 &0\end{pmatrix}
 \cdots
 \begin{pmatrix} a_j &1\\ 1 &0\end{pmatrix}
 \quad
 (j\ge 1).
\]
En effet, on v\'erifie que $|q_{j+1}|\ge \rho|q_j|$ et que $|p_{j+1}/q_{j+1}-p_j/q_j|=1/|q_jq_{j+1}|$ pour tout $j\ge 1$.  On en d\'eduit que $(p_j/q_j)_{j\ge 1}$ est une suite de Cauchy dans $K$ qui converge vers un \'el\'ement $\xi$ de $\Kbar$ avec $|\xi-p_j/q_j|\le c_0|q_j|^{-2}$ pour tout $j\ge 1$, o\`u $c_0=\rho/(\rho^2-1)$.

Comme dans \cite[\S 2]{Ro1}, on choisit $a$ et $b$ dans $A$ distincts ayant $\abs{a},\abs{b}\ge 1+\rho$, et on prend $\xi=[0,w]=[0,a,b,a,a,b,\dots]$ o\`u $w=abaab\dots$ est le mot de Fibonacci sur $\{a,b\}$, limite de la suite de mots $(w_i)_{i\ge 1}$ d\'efinie r\'ecursivement par $w_1=a$, $w_2=ab$ et $w_i=w_{i-1}w_{i-2}$ pour $i\ge 3$ dans le mono\"{\i}de $E^*$ des mots sur $\{a,b\}$.  On applique ensuite le r\'esultat de Berstel suivant lequel, pour tout $i\ge 1$, le mot $w_{i+2}$ priv\'e de ses deux derni\`eres lettres est un palindrome $m_i$.  Son image sous le morphisme de mono\"{\i}des $\Phi\colon E^*\to\GL_2(A)$ qui applique $a$ sur $\begin{pmatrix}a&1\\1&0\end{pmatrix}$ et $b$ sur $\begin{pmatrix}b&1\\1&0\end{pmatrix}$ est donc une matrice sym\'etrique qui s'\'ecrit
\[
 M_i=\Phi(m_i)=\begin{pmatrix}x_{i,0}&x_{i,1}\\x_{i,1}&x_{i,2}\end{pmatrix}
\]
pour un triplet $\ux_i=(x_{i,0},x_{i,1},x_{i,2})\in A^3$.  En vertu des consid\'erations initiales, on en d\'eduit, pour tout $i\ge 1$, les in\'egalit\'es $|x_{i,0}\xi-x_{i,1}|\le c_3|x_{i,0}|^{-1}$ et $|x_{i,1}\xi-x_{i,2}|\le c_3|x_{i,1}|^{-1}$, donc $L(\ux_i)\le c_3X_i^{-1}$ en posant $X_i=|x_{i,0}|$ et $c_3=(1+|\xi|)c_0$.

Soient $J=\begin{pmatrix}0&1\\-1&0\end{pmatrix}$ et $S = \begin{pmatrix}a&1\\1&0\end{pmatrix} \begin{pmatrix}b&1\\1&0\end{pmatrix}$.  D'apr\`es \cite[Lemme 2.1]{Ro1}, on a $M_{i+1} = M_iS_iM_{i-1}$ pour tout $i\ge 2$, o\`u $S_i=S$ si $i$ est pair et $S_i=S^t$ sinon.  Gr\^ace aux formules de \cite[\S 2]{Ro2}, on en d\'eduit que
\begin{align*}
 \det(\ux_{i-1},\ux_i,\ux_{i+1})
  &= -\trace(JM_iJM_{i+1}JM_{i-1})\\
  &=\det(M_{i-1}M_i)\det(S_iJ)
  =\pm(a-b)\neq0.
\end{align*}
Comme $x_{i,0}^{-1}\ux_i$ converge dans $\Kbar$ vers le point $(1,\xi,\xi^2)$, on conclut que $1$, $\xi$ et $\xi^2$ sont lin\'eairement ind\'ependants sur $K$ (voir la preuve de \cite[Theorem 5.1]{Ro2}).

Partant de l\`a, on reprend la preuve de \cite[Th\'eor\`eme 2.2]{Ro1}.  On en d\'eduit que les rapports $X_i/(X_{i-2}X_{i-1})$ convergent vers $|\theta|$ avec $\theta=\xi^2+(a+b)\xi+(ab+1)\neq 0$, puis qu'il existe des constantes $c_5>c_4>0$ telles que les rapports $r_i=X_i/X_{i-1}^\gamma$ satisfassent $c_4 r_{i-1}^{-1/\gamma}\le r_i\le c_5 r_{i-1}^{-1/\gamma}$ pour tout $i\ge 3$.  En choisissant $c_4$ assez petit et $c_5$ assez grand de sorte que $c_4^\gamma/c_5 \le r_2 \le c_5^\gamma/c_4$, on en d\'eduit par r\'ecurrence que $c_4^\gamma/c_5 \le r_i \le c_5^\gamma/c_4$ pour tout $i\ge 2$.  Cela corrige une erreur \`a la fin de la d\'emonstration de \cite[Th\'eor\`eme 2.2]{Ro1} et montre en particulier que $X_{i+1}\le (c_5^\gamma/c_4)X_i^\gamma$ pour tout $i\ge 1$, donc $L(\ux_i)\le c_3 X_i^{-1}\le c_2 X_{i+1}^{-1/\gamma}$ pour $i\ge 1$,
avec $c_2=c_3c_4^{1/\gamma}/c_5$.  Enfin, pour tout $X\in\bR$ avec $X\ge X_1$, il existe un indice $i\ge 1$ tel que $X_i\le X\le X_{i+1}$ et alors le point $\ux=\ux_i$ satisfait les in\'egalit\'es \eqref{ineq2}.  Donc $\xi$ poss\`ede les propri\'et\'es requises.

\medskip
\noindent
\textbf{Remerciements.} Les deux auteurs remercient le CRSNG pour son appui financier \`a cette recherche, en particulier pour la bourse de recherche de premier cycle accord\'ee au premier auteur.

\medskip
\noindent
Samuel Pilon et Damien Roy\\
D\'epartement de math\'ematiques et de statistique\\
Universit\'e d'Ottawa\\
585 King Edward, Ottawa, Ontario K1N 6N5, Canada

\smallskip
\noindent
Courriels: Samuel Pilon: \email{spilo077@uottawa.ca}, Damien Roy: \email{droy@uottawa.ca}



\begin{thebibliography}{DaSc}
%
\bibitem[Be]{Be}
  P.~Bel,
  \textit{Approximation simultan\'ee d'un nombre $v$-adique et de son carr\'e par des nombres alg\'ebriques},
  J.\ Number\ Theory {\bf 133} (2013), 3362--3380.
%
\bibitem[Bu]{Bu}
  Y.~Bugeaud,
  \textit{On simultaneous uniform approximation to a $p$-adic number and its square},
  Proc.\ Amer.\ Math.\ Soc.\ {\bf 138} (2010), 3821--3826.
%
\bibitem[DaSc]{DaSc}
  H.~Davenport, W.~M.~Schmidt,
  \textit{Approximation to real numbers by algebraic integers}, Acta Arith.\ {\bf 15}
  (1969), 393--416.
%
\bibitem[Ro1]{Ro1}
  D.~Roy,
  \textit{Approximation simultan\'ee d'un nombre et de son carr\'e},
  C.\ R.\ Acad.\ Sci.\ Paris {\bf 336} (2003), 1--6.
%
\bibitem[Ro2]{Ro2}
  D.~Roy,
  \textit{Approximation to real numbers by cubic algebraic integers I},
  Proc.\ London Math.\ Soc.\ {\bf 88} (2004), 42--62.
%
\bibitem[Te]{Te}
  O.~Teuli\'e,
  \textit{Approximation d'un nombre $p$-adique par des nombres alg\'ebriques},
  Acta Arith.\ {\bf 102} (2002), 137--155.
\bibitem[Ze]{Ze}
  D.~Zelo,
  \textit{Simultaneous approximation to real and $p$-adic numbers},
  th\`ese de doctorat,
  Universit\'e d'Ottawa, 2009; arXiv:0903.0086 [math.NT].
%
\end{thebibliography}
\end{document}